\documentclass{article}
\usepackage{amsthm}
\usepackage{amsfonts}
\usepackage{amsmath}

\newtheorem{proposition}{Proposition}
\newtheorem{theorem}{Theorem}
\newtheorem{corollary}{Corollary}
\newtheorem{lemma}{Lemma}
\newtheorem{remark}{Remark}

\def\address{Kyiv Taras Shevchenko University,
cybernetics department,Volodymyrska, 64, Kyiv, 01033, Ukraine\\
Institute of Mathematics National Academy of Sci., Tereschenkivska, 
3, Kyiv, 01601, Ukraine}

\def\email{\\prosk@imath.kiev.ua\\ Yurii\_Sam@imath.kiev.ua}

\title{Stability of the $C^*$-algebra associated with the twisted CCR.}
\author{Daniil Proskurin and Yurii Samoilenko}
\date{}
\begin {document}

\maketitle

\begin{abstract}
The universal enveloping $C^*$-algebra $\mathbf{A}_{\mu}$ of twisted canonical
commutation relations is considered. It is shown that for any
$\mu\in (-1,1)$ the $C^*$-algebra $\mathbf{A}_{\mu}$
is isomorphic to the $C^*$-algebra $\mathbf{A}_0$ generated by 
partial isometries 
$t_i,\ t_i^*,\ i=1,\ldots,d$ satisfying the relations
\[
t_i^* t_j=\delta_{ij}(1-\sum_{k<i}t_k t_k^*),\ t_j t_i=0,\ i\ne j.
\]
It is
proved that Fock representation of $\mathbf{A}_{\mu}$ is faithful.
\end{abstract}

\section*{Introduction}
Recently the interest to the *-algebras defined by generators and 
relations, their representations, particulary faithful representations, 
and  the universal 
enveloping $C^*$ -algebras has been growing because of their applications in 
mathematical physics, operator theory etc. 

A lot of interesting classes of a *-algebras depending on 
the parameters are constructed as a deformations of canonical 
commutation relations of quantum mechanics (CCR).
A well-known examples of a such deformations are
\begin{itemize}
\item $q_{ij}$-CCR introduced by M. Bozejko and R. Speicher ( see \cite{bs})
\begin{equation*}
\mathbb{C}\bigl<a_,
\ a_i^*\mid a_i^*a_j=\delta_{ij}1+q_{ij}a_ja_i^*,\ i,j=1,\ldots,d,\ 
q_{ji}=\overline{q}_{ij}\in\mathbb{C},\ \mid q_{ij}\mid\le 1\bigr>
\end{equation*}
and 
\item Twisted canonical commutation relations (TCCR) constructed 
by W. Pusz and S.L. Woronowicz ( see \cite{pw}). The TCCR have the 
following form
\begin{align}\label{mucc}
a_i^*a_i &=1+\mu^2 a_ia_i^* -(1-\mu^2)\sum_{k<i}a_k a_k^*,\
i=1,\ldots,d \nonumber \\
a_i^*a_j &=\mu a_ja_i^*,\ i\ne j,\quad
a_j a_i=\mu a_i a_j,\ i<j,\ 0<\mu<1
\end{align}
\end{itemize}

The universal $C^*$-algebra $A_{\{q_{ij}\}}$ for $q_{ij}$-CCR , 
$\mid q_{ij}\mid<\sqrt{2}-1$, was studied in ~\cite{jsw}.
 Particulary it was shown that under above restrictions on 
the coefficients $A_{\{q_{ij}\}}$ is isomorphic to the Cuntz-Toeplitz
algebra generated by isometries $\{s_i,\ s_i^*,\ i=1,\ldots,d\}$ 
satisfying relations $s_i^*s_j=0,\ i\ne j$. This implies that Fock 
representation of $A_{\{q_{ij}\}}$  is faithful. The conjecture that the 
same results are true for any choise of $\mid q_{ij}\mid<1$ was discussed 
in ~\cite{jsw} also. When
$\mid q_{ij}\mid=1,\ i\ne j$, the universal $C^*$-algebra is isomorphic
to the extension of noncommutative higher-dimesional torus 
generated by isometries satisfying $s_i^*s_j=q_{ij}s_js_i^*,\ i\ne j$ 
and the Fock representaion is faithful also (see ~\cite{dpl}).

In the present paper we consider the universal $C^*$-algebra 
$\mathbf{A}_{\mu}$ corresponding to the TCCR. Recall that the 
irreducible representations of TCCR, including unbounded, 
were described in $\cite{pw}$ and for any bounded representation 
$\pi$ of TCCR
\[
\Vert \pi(a_ia_i^*)\Vert\le\frac{1}{1-\mu^2},
\]
i.e. TCCR generate a *-bounded *-algebra ( see, for example \cite{khel} and \cite{os}).
We show in the Sec. 1 that $\mathbf{A}_{\mu}\simeq \mathbf{A}_0$ for any 
$\mu\in (-1,1)$. Note that $\mathbf{A}_0$ is generated by partial isometries 
$\{s_i,\ s_i^*,\ i=1,\ldots,d\}$ satisfying the relations 
\[
s_i^*s_j=\delta_{ij}(1-\sum_{k<i}s_k s_k^*).
\]
In the Sec.2 we prove that Fock representation of $\mathbf{A}_{\mu}$ 
is faithful.
\begin{remark}
It follows from the main result of $\cite{jps}$ that Fock representations of 
a *-algebras generated by $q_{ij}$-CCR, $\mid q_{ij}\mid<1$,
 and TCCR are faithful. In the case when $\mid q_{ij}\mid=1,\ i\ne j$, 
the kernel of Fock representation is generated as a *-ideal by the 
family $\{a_j a_i-q_{ij}a_ia_j\}$.  
\end{remark}

Finally, let us recall that by the universal $C^*$-algebra 
for a certain *-algebra $\mathcal{A}$
we mean the $C^*$-algebra $\mathbf{A}$ with the homomorphism
$\psi\colon\mathcal{A}\rightarrow\mathbf{A}$ such that for any
homomorphism $\varphi\colon\mathcal{A}\rightarrow B$, where $B$ is a
$C^*$-algebra, there exists $\theta\colon\mathbf{A}\rightarrow B$
satisfying $\theta\psi=\varphi$. It can be
obtained by the completion of $\mathcal{A}/J$ by the following 
$C^*$-seminorm on $\mathcal{A}$
\[
\Vert a\Vert=\sup_{\pi}\Vert \pi(a)\Vert,
\]
where $\sup$ is taken over all bounded representations of
$\mathcal{A}$ and $J$ is the kernel of this seminorm. Obviously this process 
requires the condition
$\sup_{\pi}\Vert \pi(a)\Vert<\infty$ for any $a\in\mathcal{A}$.

Through the paper we suppose that all $C^*$-algebras are realised by
the Hilbert space operators. Particulary, it is correct to consider
the polar decomposition of elements of a $C^*$-algebra. Obviously, we
do not claim that in general the partial isometry from the polar
decomposition lies in this $C^*$-algebra.

\section{Stability of $\mu$-CCR.}
Let us recall some properties of the $C^*$-algebra
generated by one-dimensional $q$-CCR. Namely we need the following
proposition ( see ~\cite{dpl}).
\begin{proposition} \label{qc}
Let $B$ be the unital $C^*$-algebra generated by the elements
$a,\ a^*$ satisfying the relation
\[
a^* a=1+q a a^*,\quad -1<q<1
\]
and $a^*=S^*C$ is a polar decomposition. Then $S\in B$, $B=C^*(S,S^*)$
and
\[
a=\bigl(\sum_{n=1}^{\infty} q^{n-1}S^n S^{*n}\bigr)^{\frac{1}{2}}S.
\]
\end{proposition}
Let us show that any $C^*$-algebra generated by the operators satisfying
(~\ref{mucc}) can be generated by some family of partial isometries.
\begin{proposition}
Let $A_{\mu}$   be  the  unital $C^*$-algebra  generated by  operators
$a_i,\  a_i^*,\ i=1,\ldots,d$, satisfying relations (~\ref{mucc}).
Let $a_i^*=S_i^*C_i$ be the polar decomposition.  Construct the
following family of partial isometries inductively:
\[
\widehat{S}_1:= S_1,\quad
\widehat{S}_i=(1-\sum_{j<i}\widehat{S}_j\widehat{S}_j^*)S_i
\]
Then $\forall i=1,\ldots,d$ we have $\widehat{S}_i\in A_{\mu}$,
$A_{\mu}=C^*(\widehat{S}_i,\ \widehat{S}_i^*,i=1,\ldots,d)$  and   the
following relations hold
\begin{align}\label{pisom}
\widehat{S_i}^*\widehat{S}_j &=\delta_{ij}\bigl(1-\sum_{j<i}
\widehat{S}_j\widehat{S}_j^*\bigr)\quad i,j=1,\ldots,d\\
\widehat{S} _j\widehat{S}_i & =0,\ j>i.\nonumber
\end{align}
\end{proposition}
\begin{proof}
We use induction on the number of generators.\\
{$\mathbf{d=1}$}.\\
In this case we have $a_1^*a_1=1+\mu^2 a_1a_1^*$, $a_1^*=S_1^*C_1$ and
as shown in Proposition ~\ref{qc} we have $S_1\in C^*(a_1,a_1^*)$ and
\[
a_1^*=S_1^*
\bigl(\sum_{n=1}^{\infty}\mu^{2(n-1)}S_1^n S_1^{*n}\bigr)^{\frac{1}{2}}
\]
with $S_1^*S_1=1$.\\
{$\mathbf{d-1\rightarrow d}$}.\\
Denote by $a_i^{(1)}:=(1-S_1S_1^*)a_i$,  $i=2,\ldots,d$. Note, that
the relations (~\ref{mucc}) are equivalent to
\begin{align*}
C_i^2 S_i & =S_i(1+\mu^2 C_i^2-(1-\mu^2)\sum_{j<i}C_j^2)\\
C_i^2S_j & =\mu^2 S_jC_i^2,\ j<i\\
C_i^2 S_j & = S_j C_i^2,\ j>i\\
C_iC_j &=C_jC_i,\ S_i^*S_j=S_jS_i^*,\ S_iS_j=S_jS_i
\end{align*}
Then it   is   easy  to  see  that  $(1-S_1S_1^*)a_i=a_i(1-S_1S_1^*)$,
$i=2,\ldots,d$ and
\begin{align*}
a_i^{*(1)}a_j^{(1)}& =(1-S_1S_1^*)a_i^*(1-S_1S_1^*)a_j\\
&=(1-S_1S_1^*)a_i^*a_j(1-S_1S_1^*)\\
&=\mu (1-S_1S_1^*)a_ja_i^*(1-S_1S_1^*)\\
& =\mu (1-S_1S_1^*)a_j(1-S_1S_1^*)a_i\\
&=\mu a_j^{(1)}a_i^{*(1)},\ i\ne j
\end{align*}
Analoguosly $a_j^{(1)}a_i^{(1)}=\mu a_i^{(1)}a_j^{(1)}$, $j>i>1$.
Multiplying the relation
\[
a_i^*a_i =1+\mu^2 a_ia_i^* -(1-\mu^2)\sum_{k<i}a_ka_k^*
\]
by $1-S_1S_1^*$ we get
\[
a_i^{*(1)}a_i^{(1)} =(1-S_1S_1^*)+\mu^2
a_i^{(1)}a_i^{*(1)} -(1-\mu^2)\sum_{2\le k<i}a_k^{(1)}a_k^{*(1)}.
\]
Evidently, the element $1-S_1S_1^*$ is the unit of the $C^*$-algebra  of
operators $C^*(1-S_1 S_1^*,\ a_i^{(1)},\ a_i^{*(1)},\ i=2,\ldots,d)$.
Using the assumption of induction we conclude that
\[
C^*(S_1,\ S_1^*,\ a_i^{(1)},\ a_i^{*(1)},\ i=2,\ldots,d)=
C^*(S_1,\ S_1^*,\ \widehat{S}_i^{*(1)},\
\widehat{S}_i^{(1)},\ i=2,\ldots,d)
\]
and partial isometries $\widehat{S}_i^{(1)}$, $i=2,\ldots,d$, satisfy the
relations (~\ref{pisom}). Note that
$\widehat{S}_i^{(1)}=\widehat{S_i}$, $i=2,\ldots,d$. Indeed, evidently
if $a_i^*=S_i^*C_i$ is a polar decomposition then
$a_i^{*(1)}=(1-S_1S_1^*)S_i^*(1-S_1S_1^*)C_i$, $i=2,\ldots,d$, is a polar
decomposition too. I.e. $S_i^{(1)}=(1-S_1S_1^*)S_i$
and we have
$\widehat{S}_2^{(1)}:=S_2^{(1)}=(1-S_1S_1^*)S_2=\widehat{S}_2$, further
\begin{align*}
\widehat{S}_i^{(1)}&:=(1-S_1S_1^*-\widehat{S}_2^{(1)}
\widehat{S}_2^{*(1)}
-\cdots-
\widehat{S}_{i-1}^{(1)}\widehat{S}_{i-1}^{*(1)})S_i^{(1)}\\
& = (1-S_1S_1^*-\widehat{S}_2\widehat{S}_2^*
-\cdots-
\widehat{S}_{i-1}\widehat{S}_{i-1}^*)(1-S_1S_1^*)S_{i}\\
&= (1-S_1S_1^*-\widehat{S}_2\widehat{S}_2^*
-\cdots-
\widehat{S}_{i-1}\widehat{S}_{i-1}^*)S_{i}=\widehat{S}_i\\
\end{align*}
Obviously the conclusion above is obtained by the induction.
Then $\widehat{S}_1^*\widehat{S}_i=S_1^*(1-S_1S_1^*)\widehat{S}_i=0$
and, analoguosly,
$\widehat{S}_i\widehat{S}_1=0$, $i=2,\ldots,d$. It remains only to
show that
$C^*(a_i,\ a_i^*,\ i=1,\ldots,d)=
C^*(\widehat{S}_i,\ \widehat{S}_i^*,\ i=1,\ldots,d)$. It follows from
the assumption of induction and the decomposition
\[
a_i=\sum_{n=0}^{\infty}\mu^n S_1^n a_i^{(1)}S_1^{*^n}
\]
The equality above follows from the
$S_1^*a_i=\mu a_i S_1^*$, then $S_1^{*n}a_i=\mu^n a_iS_1^{*n}$ and
\begin{align*}
&\mu^n S_1^n a_i^{(1)}S_1^{*n} =
\mu^n S_1^n(1-S_1S_1^*)a_iS_1^{*n}\\
& = S_1^n(1-S_1S_1^*)S_1^{*n} a_i
=(S_1^n S_1^{*n}-S_1^{n+1}S_1^{*n+1})a_i
\end{align*}
\end{proof}
Now we have to prove the converse statement, i.e. that any
$C^*$-algebra generated by partial isometries satisfying (~\ref{pisom})
can be generated by the elements satisfying (~\ref{mucc}).
Let us consider the unital $C^*$-algebra $A_0$ generated by
the operators $t_i,\ t_i^*,\ i=1,\ldots,d$, satisfying relations
(~\ref{pisom}). Note that $t_i,\ i=1,\ldots,d$, are partial isometries.
Indeed we have
\[
t_it_i^*t_i=t_i(1-\sum_{j<i}t_jt_j^*)=t_i.
\]
For any $i=1,\ldots,d$ define a family $\{a_i^{(j)},\ j=1,\ldots
,i\}$ inductively:
\begin{align}\label{afam}
a_i^{(i)}&=(\sum_{n=1}^{\infty}
\mu^{2(n-1)}t_i^n t_i^{*n})^{\frac{1}{2}}t_i,\\
a_i^{(j)}&=\sum_{n=0}^{\infty}\mu^n t_j^n a_i^{(j+1)}t_j^{*n},\
j=1,\ldots,i-1.\nonumber
\end{align}
We shall use the following evident decomposition also
\[
a_i^{(j)}=\sum_{n_j,\ldots,n_{i-1}=0}^{\infty}
\mu^{n_j+n_{j+1}+\cdots+n_{i-1}}t_j^{n_j}\cdots t_{i-1}^{n_{i-1}}
a_i^{(i)}t_{i-1}^{*n_{i-1}}\cdots t_j^{*n_j}
\]
Denote $a_i^{(1)}:= \widetilde{a}_i$. Our goal is to show that
$\widetilde{a}_i,\ \widetilde{a}_i^*$ satisfy the relations
(~\ref{mucc}) and
$\widehat{S}_i(\widetilde{a}_1,\ldots,\widetilde{a}_d)=t_i$,
$i=1,\ldots,d$.
To do it we prove a few auxiliary lemmas.
\begin{lemma}\label{l1}
$(1-t_1t_1^*-\cdots-t_jt_j^*)a_i^{(j)}=a_i^{(j+1)}$
\end{lemma}
\begin{proof}
In the following we denote $P_j:=1-\sum_{i\le j}t_jt_j^*$, $P_0:=1$.
It is easy
to see that $P_jt_k=0,\ k\le j$, and $P_jt_k=t_k,\ k>j$. Then
\[
P_j t_j^{n_j}\cdots t_{i-1}^{n_{i-1}}=
\left\{
\begin{array}{ccc}
0,\quad n_j\ne 0 & &\\
t_{j+1}^{n_{j+1}}\cdots t_{i-1}^{n_{i-1}},\quad n_j=0,\
\exists n_l\ne 0,\ j+1\le l\le i-1 & &\\
P_j,\quad n_l=0,\ l=j,\ldots,i-1 & &
\end{array}
\right.
\]
Then
\begin{align*}
&P_j a_i^{(j)}=
\sum_{n_j,\ldots,n_{i-1}=0}^{\infty}
\mu^{n_j+n_{j+1}+\cdots+n_{i-1}}P_jt_j^{n_j}\cdots t_{i-1}^{n_{i-1}}
a_i^{(i)}t_{i-1}^{*n_{i-1}}\cdots t_j^{*n_j}\\
&= P_j a_i^{(i)}+
\sum_{n_{j+1},\ldots,n_{i-1}=0,\sum_k n_k^2\ne 0}^{\infty}
\mu^{n_{j+1}+\cdots+n_{i-1}}t_{j+1}^{n_{j+1}}\cdots t_{i-1}^{n_{i-1}}
a_i^{(i)}t_{i-1}^{*n_{i-1}}\cdots t_{j+1}^{*n_{j+1}}\\
&=
\sum_{n_{j+1},\ldots,n_{i-1}=0}^{\infty}
\mu^{n_{j+1}+\cdots+n_{i-1}}t_{j+1}^{n_{j+1}}\cdots t_{i-1}^{n_{i-1}}
a_i^{(i)}t_{i-1}^{*n_{i-1}}\cdots t_{j+1}^{n_{j+1}}=a_i^{(j+1)}
\end{align*}
Where we have used that $P_j a_i^{(i)}=a_i^{(i)},\ j<i$. Indeed
\[
a_i^{(i)}= T_it_i,\ T_i^2=\sum_{n=1}^{\infty}\mu^{2(n-1)} t_i^nt_i^{*n}
\]
and $t_kt_k^* T_i^2=T_i^2 t_kt_k^*=0$, $i\ne k$ implies $t_k^*T_i=0$,
$i\ne k$, hence $t_k^*a_i^{(i)}=0$ and
\[
P_j a_i^{(i)}=(1-\sum_{k\le j}t_kt_k^*)a_i^{(i)}=a_i^{(i)},\quad j<i.
\]
\end{proof}
\begin{corollary}
$P_k a_i^{(j+1)}=a_i^{(j+1)},\ k\le j$
\end{corollary}
\begin{proof}
We note only that $P_kP_j=P_j,\ k\le j$.
\end{proof}
\begin{lemma}\label{l2}
$t_k^*a_i^{(j+1)}=0,\ a_i^{(j+1)}t_k=0$,
$t_k^*a_i^{*(j+1)}=0,\ a_i^{*(j+1)}t_k=0$, for any $k\le j<i$.
\end{lemma}
\begin{proof}
As in the previous lemma we have $t_k^*a_i^{(j+1)}=t_k^*a_i^{(i)}=0$
and $a_i^{(j+1)}t_k=a_i^{(i)}t_k=T_it_it_k=0$ since $t_it_k=0,\ i>k$.
The other relations are adjoint to the proved above.
\end{proof}
\begin{lemma}\label{l3}
\[
t_j^{*n}t_j^m=\left\{
\begin{array}{ccc}
t_j^{*n-m},\ n>m &&\\
P_{j-1},\ n=m &&\\
t_j^{m-n},\  n<m &&
\end{array}
\right.
\]
\end{lemma}
\begin{proof}
Induction on $n,m$ using the basic relations (~\ref{pisom}).
\end{proof}
Now we are able to prove the following proposition.
\begin{proposition}\label{prop1}
For any $i=1,\ldots,d$ and $1\le j\le i$ we have
\[
a_i^{*(j)}a_i^{(j)}=P_{j-1}+\mu^2a_i^{(j)}a_i^{*(j)}-
(1-\mu^2)\sum_{j\le k<i}a_k^{(j)}a_k^{*(j)}
\]
\end{proposition}
\begin{proof}
We use the induction on $j$ for a fixed $i=1,\ldots,d$.\\
For $\mathbf{j=i}$ we have
\begin{align*}
& a_i^{*(i)}a_i^{(i)}=t_i^*\bigl(\sum_{n=1}^{\infty}
\mu^{2(n-1)}t_i^nt_i^{*n}\bigr)t_i\\
&= t_i^*t_i+\mu^2\sum_{n=1}^{\infty}\mu^{2(n-1)}
t_i^nt_i^{*n}\\
&= P_{i-1}+\mu^2 a_i^{(i)}a_i^{*(i)}
\end{align*}
$\mathbf{j+1\rightarrow j}$.
Using the results of previous lemmas we have
\begin{align*}
& a_i^{*(j)}a_i^{(j)}=\sum_{n,m=0}^{\infty}\mu^{n+m}
t_j^n a_i^{*(j+1)}t_j^{*n}t_j^m a_i^{(j+1)}t_j^{*m}
= \sum_{n=0}^{\infty}\mu^{2n}
t_j^n a_i^{*(j+1)} a_i^{(j+1)}t_j^{*n}\\
&= \sum_{n=0}^{\infty}\mu^{2n}
t_j^n\bigl( P_j+\mu^2 a_i^{(j+1)} a_i^{*(j+1)}
-(1-\mu^2)\sum_{j+1\le k<i}a_k^{(j+1)}a_k^{*(j+1)}\bigr)t_j^{*n}\\
&=\sum_{n=0}^{\infty}\mu^{2n}t_j^n P_j t_j^{*n}+
\mu^2 a_i^{(j)}a_i^{*(j)}-
(1-\mu^2)\sum_{j+1\le k<i}a_k^{(j)}a_k^{*(j)}\\
&=P_j+\sum_{n=1}^{\infty}\mu^{2n}(t_j^n t_j^{*n}-t_j^{n+1}t_j^{*n+1})
+\mu^2 a_i^{(j)}a_i^{*(j)}-
(1-\mu^2)\sum_{j+1\le k<i}a_k^{(j)}a_k^{*(j)}\\
&=P_{j-1}-(1-\mu^2)\sum_{n=1}^{\infty}\mu^{2(n-1)}t_j^n t_j^{*n}+
\mu^2 a_i^{(j)}a_i^{*(j)}-
(1-\mu^2)\sum_{j+1\le k<i}a_k^{(j)}a_k^{*(j)}\\
&=P_{j-1}+\mu^2 a_i^{(j)}a_i^{*(j)}-
(1-\mu^2)\sum_{j\le k<i}a_k^{(j)}a_k^{*(j)}
\end{align*}
Particulary, for $j=1$ we have
\[
\widetilde{a}_i^* \widetilde{a}_i=
1+\mu^2 \widetilde{a}_i \widetilde{a}_i^*-
(1-\mu^2)\sum_{k<i}\widetilde{a}_k \widetilde{a}_k^*
\]
\end{proof}
It remains to show that
$\widetilde{a}_i^*\widetilde{a}_j=\mu \widetilde{a}_j\widetilde{a}_i$.
Then
$\widetilde{a}_j\widetilde{a}_i=\mu \widetilde{a}_i\widetilde{a}_j$, $j>i$,
hold automatically (see ~\cite{jsw1}).
\begin{lemma}
$a_i^{*(k)}a_j^{(j)}=0$, $j< k\le i$.
\end{lemma}
\begin{proof}
We use induction again. For $\mathbf{k=i}$ one has
\[
a_i^{*(i)}a_j^{(j)}=t_i^* T_i T_j t_j=0
\]
since $T_i^2 T_j^2 =0,\ i\ne j$, and $T_i,\ T_j\ge 0$.\\
$\mathbf{k+1\rightarrow k}$.
\[
a_i^{*(k)}a_j^{(j)}=\sum_{n=0}^{\infty}\mu^n
t_k^n a_i^{*(k+1)}t_k^{*n} a_j^{(j)}
= a_i^{*(k+1)}a_j^{(j)}=0
\]
since $t_k^* a_j^{(j)}=0$, $k>j$ (see Lemma ~\ref{l1}).
\end{proof}
\begin{lemma}
$a_i^{*(j)}a_j^{(j)}=\mu a_j^{(j)}a_i^{*(j)}$, $j<i$.
\end{lemma}
\begin{proof}
Let us show that $a_i^{*(j)}T_j^2=T_j^2 a_i^{*(j)}$ and
$a_i^{*(j)}t_j=\mu t_j a_i^{*(j)}$.
\begin{align*}
& a_i^{*(j)}T_j^2=\bigl(\sum_{n=0}^{\infty}
\mu^n t_j^n a_i^{*(j+1)} t_j^{*n}
\bigr)
\bigl(
\sum_{m=1}^{\infty}\mu^{2(m-1)}t_j^m t_j^{*m}
\bigr)\\
&=\sum_{n=0,m=1}^{\infty}\mu^{n+2(m-1)}t_j^n
a_i^{*(j+1)}t_j^{*n}t_j^m t_j^{*m}\\
&= \sum_{n=1,m\le n}^{\infty}\mu^{n+2(m-1)}
t_j^n a_i^{*(j+1)}t_j^{*n}.
\end{align*}
where we have used Lemmas ~\ref{l2},\ref{l3}. Analogously
\[
T_j^2 a_i^{*(j)}=
\sum_{n=1,m\le n}^{\infty}\mu^{n+2(m-1)}
t_j^n a_i^{*(j+1)}t_j^{*n}= a_i^{*(j)}T_j^2.
\]
Finally
\begin{align*}
& a_i^{*(j)}t_j=\bigl(
\sum_{n=0}^{\infty}\mu^n t_j^n a_i^{*(j+1)}t_j^{*n}
\bigr)t_j\\
&=\mu t_j a_i^{*(j+1)}t_j^*t_j+
\sum_{n=2}^{\infty}\mu^n t_j^n a_i^{*(j+1)}t_j^{*n-1}\\
&=\mu t_j a_i^{*(j+1)}P_{j-1}+\mu t_j
\sum_{n=1}^{\infty}\mu^n t_j^n a_i^{*(j+1)}t_j^{*n}\\
&=\mu t_j\sum_{n=0}^{\infty}\mu^n t_j^n a_i^{*(j+1)}t_j^{*n}
=\mu t_j a_i^{*(j)}.
\end{align*}
Then
\[
a_i^{*(j)}a_j^{(j)}=a_i^{*(j)}T_jt_j=\mu T_jt_ja_i^{*(j)}
=\mu a_j^{(j)}a_i^{*(j)},\ i>j.
\]
\end{proof}
\begin{lemma}
$a_i^{*(k)}a_j^{(k)}=\mu a_j^{(k)}a_i^{*(k)}$, $1\le k<j<i$
\end{lemma}
\begin{proof}
We use induction. The case $k=j$ is considered in the Lemma above.\\
$\mathbf{k+1\rightarrow k}$. As in the Proposition ~\ref{prop1} we have
\begin{align*}
& a_i^{*(k)}a_j^{(k)}=\sum_{n=0}^{\infty}
\mu^{2n} t_k^n a_i^{*(k+1)}a_j^{(k+1)}t_k^{*n}\\
& =\mu \sum_{n=0}^{\infty}
\mu^{2n} t_k^n a_j^{(k+1)} a_i^{*(k+1)}t_k^{*n}
=\mu a_j^{(k)}a_i^{(k)}.
\end{align*}
Particulary, for $k=1$ we have
$\widetilde{a}_i^* \widetilde{a}_j=\mu \widetilde{a}_j\widetilde{a}_i^*$,
$i>j$.
\end{proof}
So, we have proved the following theorem.
\begin{theorem}
Let $A_0=C^*(t_i,\ t_i^*,\ i=1,\ldots,d)$ where
$\{t_i,\ t_i^*,\ i=1,\ldots,d\}$
satisfy relations (~\ref{pisom}), and the family
i$\{\widetilde{a}_i, \widetilde{a}_i^*,\ i\ge 1\}$
is constructed according to formulas
(~\ref{afam}). Then the relations (~\ref{mucc}) are satisfied and
we have
$\widehat{S}_i(\widetilde{a}_1,\ldots,\widetilde{a}_d)=t_i,\ i=1,\ldots,d$.
\end{theorem}
\begin{corollary}
For any $\mu\in (-1,1)$ the $C^*$-algebra $\mathbf{A}_{\mu}$ is isomorphic
to $\mathbf{A}_0$.
\end{corollary}
\begin{proof}
Using the universal property of $\mathbf{A}_0$ we can define the
surjective homomorphism
$\varphi\colon \mathbf{A}_0\rightarrow \mathbf{A}_{\mu}$ by rule
$\varphi(t_i)=\widehat{S_i}$, $i=1,\ldots,d$. Analogously, we have
$\psi\colon \mathbf{A}_{\mu}\rightarrow \mathbf{A}_0$,
$\psi (a_i)= a_i^{(1)}$,
$i=1,\ldots,d$. Obviously, $\psi\varphi=id$ and $\varphi\psi=id$.
\end{proof}
\section{Fock representation.}
Recall that Fock representation of TCCR is the irreducible
representation determined by the cyclic vector $\Omega$ such that
$a_i^*\Omega=0$, $i=1,\ldots,d$.

Let us prove that Fock representation of $\mathbf{A}_{\mu}$ is faithful.
Firstly note that Fock representation of $\mathbf{A}_0$ corresponds to
the Fock representation of $\mathbf{A}_{\mu}$
(it can be easely seen from the formulas
connecting $\{t_j\}$ and $\{a_j\}$).

In the following we need the description of classes of unitary
equivalence of irreducible representations of $\mathbf{A}_0$. As we
have noted above, the
irreducible representations of TCCR,
including unbounded representations, were classified in ~\cite{pw}.
However it is more
convenient for us to present the representations of $\mathbf{A}_0$ in
some different form .
\begin{proposition}
Let $\pi$ be an irreducible representation of $\mathbf{A}_0$ acting on the
Hilbert space $\mathcal{H}$, then for some
$j=1,\ldots,d$ we have
$\mathcal{H}\simeq \bigotimes_{k=1}^j l_2(\mathbb{N})$ and
\begin{align*}
\pi(t_i)& =\bigotimes_{k=1}^{i-1}(1-S S^*)\otimes S\otimes
\bigotimes_{k=i+1}^j 1,\ i\le j\\
\pi(t_{j+1})& = e^{i\varphi}\bigotimes_{k=1}^j (1-S S^*),\
\varphi\in [0,2\pi)\\
\pi(t_i)&=0,\ i>j+1,
\end{align*}
where $S$ is a unilateral shift on $l_2(\mathbb{N})$.
The case $j=d$ corresponds to the Fock representation.
\end{proposition}
\begin{proof}
It follows from (\ref{pisom}) that $\pi(t_1)$ is isometry. Hence,
either $\ker\pi(t_1^*)\ne\{0\}$ or $\pi(t_1)$ is unitary. We shall use
here $t_i$ instead $\pi(t_i)$.

Let $\ker t_1^*=\mathcal{H}_1\ne\{0\}$.
Then the relations (~\ref{pisom}) imply that
$\mathcal{H}=\bigoplus_{n=0}^{\infty} t_1^n\mathcal{H}_1$ and
\[
t_i,\ t_i^*\colon\mathcal{H}_1\rightarrow\mathcal{H}_1,\
t_i,\ t_i^*\colon t_1^{n}\mathcal{H}_1\rightarrow \{0\},\ n>1,\ i>1.
\]
If we identify $t_1^n\mathcal{H}_1$ with $e_n\otimes\mathcal{H}_1$,
$n\ge 0$, then
\[
t_1=S\otimes 1,\ t_i=(1-S S^*)\otimes t_i^{(1)},\ i>1
\]
where $S e_n=e_{n+1},\ n\ge 0$ and the family $\{t_i^{(1)},\ i>1\}$
saisfy (~\ref{pisom}) on the space $\mathcal{H}_1$. Moreover, it is
easy to show that the family $\{t_i,\ i=1,\ldots,d\}$ is irreducible iff
$\{t_i^{(1)},\ i=2,\ldots,d\}$ is irreducible.

If $t_1$ is unitary, then $t_it_1=0$, $i>1$, implies $t_i=0$.
\end{proof}
Using the previous proposition we can prove the following theorem.
\begin{theorem}
The Fock representation of $\mathbf{A}_0$ is faithful.
\end{theorem}
\begin{proof}
Let $C_F$ be the $C^*$-algebra generated by the operators of Fock
representation and $C_{\pi}$ be the $C^*$-algebra generated by some
irreducible representation $\pi$ of $\mathbf{A}_0$.
To prove the statement it is sufficient to construct a
homomorphism
\[
\psi\colon C_F\rightarrow C_{\pi}
\]
such that $\pi=\psi\pi_F$ (then $\pi(\ker\ \pi_F)=\{0\}$ for any
irreducible representation of $\mathbf{A}_0$, i.e. $\ker\ \pi_F=\{0\}$, where
we denote by $\pi_F$ the Fock representation).

To do it, we note that if $\pi$ corresponds to some $j=1,\ldots,d-1$,
then $C_F$ and $C_{\pi}$ are the $C^*$-subalgebras of the
$\bigotimes_{k=1}^d C^*(S,S^*)$ and $\bigotimes_{k=1}^j C^*(S,S^*)$
respectively.
Recall that $C^*(S,S^*)\simeq\mathcal{T}(C(\mathbf{T}))$ is a
nuclear $C^*$-algebra of the Toeplitz
operators. Then we can define the homomomorphism
\[
\psi\colon
\bigotimes_{k=1}^d C^*(S,S^*)\rightarrow\bigotimes_{k=1}^j C^*(S,S^*),
\]
defined by
\begin{align*}
\psi(\otimes_{k=1}^{i-1}1\otimes S\otimes_{k=i+1}^d 1)&=
\otimes_{k=1}^{i-1}1\otimes S\otimes_{k=i+1}^j 1,\ i\le j\\
\psi(\otimes_{k=1}^{j}1\otimes S\otimes_{k=i+1}^d 1)&=
e^{i\varphi}\otimes_{k=1}^{j}1,\\
\psi(\otimes_{k=1}^{i-1}1\otimes S\otimes_{k=i+1}^d 1)&=
\otimes_{k=1}^{j}1,\ i> j\\
\end{align*}
It remains only to restrict $\psi$ onto $C_F$ and to note that
$\psi(C_F)=C_{\pi}$.
\end{proof}
{\bf{Acknoledgements.}}\\
We express our gratitude to Vasyl Ostrovsky\u{\i} and Stanislav 
Popovich for their critical remarks and helpful discussions.

This work was partially supported by the State Fund of Fundamental
Researches of Ukraine, grant no. 01.07/071.

\nopagebreak
\vskip 0.5cm
\noindent
\address
\nopagebreak

\bigskip\noindent{\em e-mail: }\email


\begin{thebibliography}{99}
\bibitem{dn} { K. Dykema and A.Nica}. On the Fock representation
of the q-commutation relations. {\em J. Reine Angew. Math.} 
\textbf{440} (1993), 201--212.

\bibitem{bs}{ Bo\.zejko, M. and  Speicher, R.}: Completely positive maps on
Coxeter groups, deformed commutation relations, and operator
spaces. {\it Mat. Ann.} {\bf 300} (1994), 97--120.

\bibitem{jps} { J{\o}rgensen, P.E.T,  Proskurin, D.P. and
 Samo\u\i{}lenko, Yu.S.}: The kernel of Fock representation of Wick
algebras with braided operator of coefficients. Accepted to
publication in {\it Pacific. Journ. Math.}, math-ph/0001011.

\bibitem{jsw} { J{\o}rgensen, P.E.T., Schmitt L.M., and Werner, R.F.}:
q-canonical commutation relations and stability of the Cuntz algebra.
{\it Pacific Journ. Math.} {\bf 163}, no. 1 (1994), 131--151

\bibitem{jsw1} {  J{\o}rgensen, P.E.T. Schmitt,L.M. and  Werner,R.F.}.
Positive representations of general commutation relations allowing Wick
ordering. {\it J.\ Funct.\ Anal.} {\bf 134} (1995), 33-99.

\bibitem{khel} {Khelemski\u\i{}, A.Ya.}: Banach algebras and poly-normed algebras:
general theory, representations, homologies. {\it Nauka, Moscow} (1989), (Russian).

\bibitem{os} { Ostrovsky\u\i{}, V. and  Samo\u\i{}lenko, Yu.}: Introduction
to the Theory of Representations of Finitely Presented *-Algebras. I.
Representations by bounded operators. {\it The Gordon and Breach
Publishing group, London} (1999).

\bibitem{pw} Pusz, W. and Woronowicz, S.L.: {Twisted second
quantization.} {\it Reports Math. Phys.}\textbf{27} (1989), 251--263.

\bibitem{dpl} Proskurin, D.: Stability of a special class of
$q_{ij}$-CCR and extensions of higher-dimensional noncommutative tori,
to be printed in {\it Lett. Math. Phys.}
\end{thebibliography}
\end{document}